\documentclass[11pt]{amsart}%
\usepackage{amsfonts}
\usepackage{amsmath}
\usepackage{amssymb}
\usepackage{graphicx}%
\providecommand{\U}[1]{\protect\rule{.1in}{.1in}}

\theoremstyle{plain}

\numberwithin{equation}{section}
\begin{document}
\title[Quadratic Hessian equation]{Quadratic Hessian equation}
\author{Yu YUAN}
\address{University of Washington\\
Department of Mathematics, Box 354350\\
Seattle, WA 98195}
\email{yuan@math.washington.edu}
\thanks{This work is partially supported by an NSF grant.}
\dedicatory{In memory of my teacher, Ding Wei-Yue laoshi (1945--2014)}\date{\today }
\maketitle

\section{Introduction}

\label{sec:Intro}

The elementary symmetric quadratic or $\sigma_{2}$ equation%
\begin{equation}
\sigma_{2}\left(  D^{2}u\right)  =\sum_{1\leq i_{1}<i_{2}\leq n}\lambda
_{i_{1}}\lambda_{i_{2}}=\frac{\left(  \bigtriangleup u\right)  ^{2}-\left\vert
D^{2}u\right\vert ^{2}}{2}=1 \label{Esigma2}%
\end{equation}
with $\lambda_{i}^{\prime}s$ being the eigenvalues of the Hessian $D^{2}u$ of
scalar function $u,$ is a nonlinear Hessian dependence equation of the lowest
integer order, and is called fully nonlinear equation, because the
nonlinearity is on the highest order derivatives of the solutions. The
$\sigma_{2}$ equation sits in between the (linear) Laplace equation
$\sigma_{1}\left(  D^{2}u\right)  =\lambda_{1}+\cdots+\lambda_{n}%
=\bigtriangleup u=1$ and the (fully nonlinear) Monge-Amp\`{e}re equation
$\sigma_{n}\left(  D^{2}u\right)  =\lambda_{1}\lambda_{2}\cdots\lambda
_{n}=\det D^{2}u=1.$ The 2-sheet hyperboloid level set of the equation%
\[
\left\{  \lambda\in\mathbb{R}^{n}:\ \lambda_{1}+\cdots+\lambda_{n}\ =\pm
\sqrt{2+\left\vert \lambda\right\vert ^{2}}\right\}
\]
is rotationally symmetric, unlike all the other $\sigma_{k}$ equations with
$3\leq k\leq n.$

To make those equations elliptic, or monotone dependence on Hessian along
positive definite symmetric matrices, we require the linearized operator
positive definite. Equivalently, the normal of the level set has positive sign
for all components in the eigenvalue space, respectively positive definite in
the matrix space. For example, among all four branches of level set
$\lambda_{1}\lambda_{2}\lambda_{3}=1,$ only one is elliptic; the same is true
for $\lambda_{1}\lambda_{2}\lambda_{3}=-1.$ The negative definite or all
negative component case is also considered as elliptic. In particular, the two
symmetric branches of $\sigma_{2}\left(  \lambda\right)  =1$ are both
elliptic. The choice of branch is made automatically by $C^{2}$ solutions. For
less smooth solutions such as continuous viscosity ones, the positive
$\bigtriangleup u>0$ or negative $\bigtriangleup u<0$ branch has to be
specified. Afterward, $\mathbf{-}D^{2}u$ is on the other branch in viscosity sense.

Replacing the flat eigenvalues $\lambda_{i}^{\prime}s$ with the principal
curvatures $\kappa_{i}^{\prime}s$ of graph $\left(  x,u\left(  x\right)
\right)  $ in Euclid space $\mathbb{R}^{n}\times\mathbb{R}^{1},$ one has the
scalar curvature equation%
\[
\sigma_{2}\left(  \kappa_{1},\cdots,\kappa_{n}\right)  =1.
\]
Recall $\kappa_{i}^{\prime}s$ are the eigenvalues of the normalized second
fundamental form $II$ by the induced metric $g$ or shape matrix%
\[
IIg^{-1}=\frac{D^{2}u}{\sqrt{1+\left\vert Du\right\vert ^{2}}}\left[
I-\frac{Du\otimes Du}{1+\left\vert Du\right\vert ^{2}}\right]  =\left[
\partial_{x_{i}}A_{p_{j}}\left(  Du\right)  \right]  ,
\]
where $II=D^{2}u/\sqrt{1+\left\vert Du\right\vert ^{2}},\ g=I+Du\otimes
Du,\ $and$\ A\left(  p\right)  =\sqrt{1+\left\vert p\right\vert ^{2}}.$
Replacing the flat eigenvalues $\lambda_{i}^{\prime}s$ with the eigenvalues of
Schouten tensor of a conformal metric $g=u^{-2}g_{0},$ one has $\sigma_{2}%
$-type Yamabe equation in conformal geometry, which simplifies to%
\[
\sigma_{2}\left(  uD^{2}u-\frac{1}{2}\left\vert Du\right\vert ^{2}I\right)
=1
\]
for flat metric $g_{0}.$ Replacing the flat eigenvalues $\lambda_{i}^{\prime
}s$ with the eigenvalues of Hermitian Hessian $\partial\bar{\partial}u,$ one
has $\sigma_{2}$-type equation arising from complex geometry. Lastly, in three
dimensions $\sigma_{2}\left(  D^{2}u\right)  =1$ or equivalently
$\arctan\lambda_{1}+\cdots+\arctan\lambda_{3}=\pm\pi/2$ is the potential
equation of minimal Lagrangian graph $\left(  x,Du\left(  x\right)  \right)  $
with phase $\pm\pi/2$ in Euclid space $\mathbb{R}^{3}\times\mathbb{R}^{3}.$

\section{Results}

\subsection{Outline}

Once an equation is given, the first question to answer is the existence of
solutions. Smooth ones cannot be obtained immediately, in general; worse, they
may not even exist. The typical approach is to first seek weak solutions, in
the integral sense if the equation has divergence structure, or in the
\textquotedblleft pointwise integration by parts sense\textquotedblright,
namely, in the viscosity sense if the equation enjoys a comparison principle.
After obtaining those weak solutions, one studies the regularity and other
properties of the solutions, such as Liouville or Bernstein type rigidity for
entire solutions. All these hinge on \textit{a priori} estimates of
derivatives of solutions:
\[
\Vert D^{2}u\Vert_{L^{\infty}(B_{1})}\leq C(\Vert Du\Vert_{L^{\infty}(B_{2}%
)})\leq C(\Vert u\Vert_{L^{\infty}(B_{3})}).
\]
Having the $L^{\infty}$ bound of the Hessian available, the ellipticity of the
above fully nonlinear equations becomes uniform, we can apply the
Evans-Krylov-Safonov theory (for the ones with convexity/concavity, possibly
without divergence structure) or the Evans-Krylov-De Giorgi-Nash theory (for
the ones with convexity/concavity and divergence structure) to obtain
$C^{2,\alpha}$ estimates of solutions. Either theory can handle the quadratic
Hessian equation along with all other $\sigma_{k}$ equations, because they all
share the divergence structure. In the $\sigma_{2}$ equation (\ref{Esigma2})
case, the linearized operator is readily seen in divergence form
\[
\bigtriangleup_{F}=\sum_{i,j=1}^{n}F_{ij}\partial_{ij}=\sum_{i,j=1}%
^{n}\partial_{i}\left(  F_{ij}\partial_{j}\right)  ,\ \
\]
with%
\[
\left(  F_{ij}\right)  =\bigtriangleup u\ I-D^{2}u=\sqrt{2+\left\vert
D^{2}u\right\vert ^{2}}\ I-D^{2}u>0.
\]
Here and in the remaining, for certainty, and without loss of generality, we
assume $\bigtriangleup u>0.$ The concavity of the equation is evident in an
equivalent form
\begin{equation}
\bigtriangleup u-\sqrt{2+\left\vert D^{2}u\right\vert ^{2}}=0, \label{Eec}%
\end{equation}
whose linearized operator%
\[
I-\frac{D^{2}u}{\sqrt{2+\left\vert D^{2}u\right\vert ^{2}}}%
\]
also reveals the uniform ellipticity of the equation with uniformly bounded
Hessian solutions.

Considering the minimal surface structure of the quadratic Hessian equation in
three dimensions, the $C^{2,\alpha}$ estimate can also be achieved via
geometric measure theory. For the $\sigma_{n}$ or Monge-Amp\`{e}re equation,
back in the 1950s, Calabi attained $C^{3}$ estimates by interpreting the cubic
derivatives in terms of the curvature of the corresponding Hessian metric
$g=D^{2}u.$ Further, iterating the classic Schauder estimates, one obtains
smoothness of the solutions, and even analyticity, if the smooth equations
such as all the $\sigma_{k}$ equations are also analytic.

\subsection{Rigidity of entire solutions}

The classic Liouville type theorem asserts every entire harmonic function
bounded from below or above is a constant by the Harnack inequality. Thus
every semiconvex harmonic function is a quadratic one, as its double
derivatives are all harmonic with lower bounds, hence constants. Similarly,
every entire (convex) solution to the Monge-Amp\`{e}re equation $\det
D^{2}u=1$ is quadratic. This was first proved in two dimensional case by
J\"{o}rgens, later in dimension up to five by Calabi, and in all dimensions by
Pogorelov. Also, Cheng-Yau had a geometric proof.

Recently, Shankar-Yuan \cite{SY3} proved that every entire semiconvex solution
to the quadratic Hessian equation $\sigma_{2}\left(  D^{2}u\right)  =1$ is
quadratic. In dimension two, it is the above classic J\"{o}rgens's theorem
without any extra condition (not even convexity) on the entire solutions, thus
a Bernstein type result. In three dimensions, this was proved in \cite{Y}
earlier, as a by-product of rigidity for the special Lagrangian equation.

Under an almost convexity condition on entire solutions to $\sigma_{2}\left(
D^{2}u\right)  =1$ in general dimension, Chang-Yuan derived the rigidity
\cite{ChY}. Under a general semiconvexity and an additional quadratic growth
assumption on entire solutions in general dimension, Shankar-Yuan showed the
rigidity in \cite{SY1}. Assuming only quadratic growth on entire solutions to
$\sigma_{2}\left(  D^{2}u\right)  =1$ in three and four dimensions, the same
rigidity result was proved in the joint work with Warren \cite{WY} and Shankar
\cite{SY4} respectively. Assuming a super quadratic growth condition,
Bao-Chen-Ji-Guan \cite{BCGJ} demonstrated that all convex entire solutions to
$\sigma_{2}\left(  D^{2}u\right)  =1\ $along with other $\sigma_{k}\left(
D^{2}u\right)  =1$ are quadratic polynomials; and Chen-Xiang \cite{CX}[ showed
that all \textquotedblleft super quadratic\textquotedblright\ entire solutions
to $\sigma_{2}\left(  D^{2}u\right)  =1$ with $\sigma_{1}\left(
D^{2}u\right)  >0$ and $\sigma_{3}\left(  D^{2}u\right)  \geq-K$ are also
quadratic polynomials.

Warren's rare saddle entire solution $u\left(  x_{1},\cdots,x_{n}\right)
=\left(  x_{1}^{2}+x_{2}^{2}-1\right)  e^{x_{3}}+\frac{1}{4}e^{-x_{3}}$ to
$\sigma_{2}\left(  D^{2}u\right)  =1$ in dimension three and above \cite{W},
confirms the necessity of the semiconvexity or the quadratic growth
assumption. C.-Y. Li \cite{L} followed with \textquotedblleft
non-degenerate\textquotedblright\ entire solution $u\left(  x\right)  =\left(
x_{1}^{2}+x_{2}^{2}-1\right)  e^{x_{n}}+\frac{n-2}{4}e^{-x_{n}}+\left(
x_{3}+\cdots+x_{n-1}\right)  x_{n}$ in dimension $n$ and above for $n\geq4.$

\subsubsection{Two dimensions}

In the following, we recap Nitsche's idea in showing the rigidity of entire
solutions in two dimensions.

Given a $C^{2}$ solution $u$ to $\sigma_{2}\left(  D^{2}u\right)  =1,$ up to
negation, we assume $D^{2}u$ is on the positive branch of the hyperbola
$\lambda_{1}\lambda_{2}=1,$ in turn, $u$ is convex.\ Let $w$ be the
Legendre-Lewy transform of $u\left(  x\right)  $, that is, the Legendre
transform of $u\left(  x\right)  +\left\vert x\right\vert ^{2}/2.\ $%
Geometrically their \textquotedblleft gradient\textquotedblright\ graphs
satisfy $\left(  x,Du\left(  x\right)  +x\right)  =\left(  Dw\left(  y\right)
,y\right)  \in\mathbb{R}^{2}\times\mathbb{R}^{2},$ and the \textquotedblleft
slopes\textquotedblright\ of graphs satisfy%
\[
\left(  I,D^{2}u\left(  x\right)  +I\right)  =\left(  D^{2}w\left(  y\right)
\frac{\partial y}{\partial x},\frac{\partial y}{\partial x}\right)  .
\]
It follows that%
\[
I<D^{2}u\left(  x\right)  +I=\left(  D^{2}w\left(  y\right)  \right)
^{-1}\ \ \text{or \ \ }D^{2}u\left(  x\right)  =\left(  D^{2}w\left(
y\right)  \right)  ^{-1}-I.
\]
Taking determinants yields%
\[
1=\sigma_{2}\left(  D^{2}u\right)  =\det\left[  \left(  D^{2}w\left(
y\right)  \right)  ^{-1}-I\right]  =\mu_{1}^{-1}\mu_{2}^{-1}-\mu_{1}^{-1}%
-\mu_{2}^{-1}+1.
\]
or an equation for the eigenvalues $\mu_{i}^{\prime}s$ of $D^{2}w$%
\[
1=\mu_{1}+\mu_{2}=\bigtriangleup w.
\]

Noticing the boundedness of Hessian $D^{2}w$%
\[
0<D^{2}w<I,
\]
we see the constancy by Liouville. Consequently from the flatness of graphs
$\left(  x,Du\left(  x\right)  +x\right)  =\left(  Dw\left(  y\right)
,y\right)  $ or constancy of $\left(  D^{2}w\right)  ^{-1}-I=D^{2}u,$ it
follows that $u$ is quadratic. Note that J\"{o}rgens' original
\textquotedblleft involved\textquotedblright\ proof made use of a partial
Legendre transformation.

Going further, this Legendre-Lewy transformation proof of J\"{o}rgens'
theorem, coupled with Heinz transformation, led Nitsche \cite{N} to his
elementary proof of the original Bernstein theorem: every entire solution to
the minimal surface equation%
\[
\operatorname{div}\left(  \frac{Df}{\sqrt{1+\left\vert Df\right\vert ^{2}}%
}\right)  =0
\]
in two dimensions is linear.

The remaining proof goes as follows. The mean curvature vector $\vec{H}$ of
the graph $\left(  x_{1},x_{2},f\left(  x\right)  \right)  $ is%
\begin{align*}
\vec{H} &  =\bigtriangleup_{g}\left(  x_{1},x_{2},f\left(  x\right)  \right)
=\sum_{i,j}\frac{1}{\sqrt{g}}\partial_{i}\left(  \sqrt{g}g^{ij}\partial
_{j}\right)  \left(  x_{1},x_{2},f\left(  x\right)  \right)  \\
&  =\left(  \operatorname{div}\frac{\left(  1+f_{2}^{2},-f_{1}f_{2}\right)
}{\sqrt{1+\left\vert Df\right\vert ^{2}}},\operatorname{div}\frac{\left(
-f_{1}f_{2},1+f_{1}^{2}\right)  }{\sqrt{1+\left\vert Df\right\vert ^{2}}%
},\operatorname{div}\frac{\left(  f_{1},f_{2}\right)  }{\sqrt{1+\left\vert
Df\right\vert ^{2}}}\right)  ,
\end{align*}
where $\sqrt{g}g^{-1}=\left[
\begin{array}
[c]{cc}%
1+f_{2}^{2} & -f_{1}f_{2}\\
-f_{2}f_{1} & 1+f_{1}^{2}%
\end{array}
\right]  /\sqrt{1+\left\vert Df\right\vert ^{2}}$ was used in the last
equality. The magnitude $H$ of the mean curvature vector is, as noted before%
\[
H=Tr\left[  \partial_{x_{i}}A_{p_{j}}\left(  Df\right)  \right]
=\operatorname{div}\frac{\left(  f_{1},f_{2}\right)  }{\sqrt{1+\left\vert
Df\right\vert ^{2}}}=0.
\]
Hence, $\bigtriangleup_{g}\left(  x_{1},x_{2},f\left(  x\right)  \right)
=\left(  0,0,0\right)  .$ Certainly $\bigtriangleup_{g}x_{1}=0=\bigtriangleup
_{g}x_{2}$ can be verified by direct computation.

Considering the first component equation $\bigtriangleup_{g}x_{1}=0,$ the
conjugate function of $x_{1}$ is defined as%
\[
x_{1}^{\ast}\left(  x_{1},x_{2}\right)  =\int^{\left(  x_{1},x_{2}\right)
}\frac{f_{1}f_{2}dx_{1}+\left(  1+f_{2}^{2}\right)  dx_{2}}{\sqrt{1+\left\vert
Df\right\vert ^{2}}}.
\]
Similarly, the conjugate function of $x_{2}$ is also defined. Together, they
represent the \textquotedblleft normalized\textquotedblright\ metric%
\[
\frac{1}{\sqrt{1+\left\vert Df\right\vert ^{2}}}\left[
\begin{array}
[c]{cc}%
1+f_{1}^{2} & f_{1}f_{2}\\
f_{2}f_{1} & 1+f_{2}^{2}%
\end{array}
\right]  =\left[
\begin{array}
[c]{c}%
Dx_{2}^{\ast}\\
Dx_{1}^{\ast}%
\end{array}
\right]  .
\]
By symmetry of the left-side matrix, $\partial_{2}x_{2}^{\ast}=\partial
_{1}x_{1}^{\ast},$ then there exists a double\ potential $u$ so that
$Du=\left(  x_{2}^{\ast},x_{1}^{\ast}\right)  .$ Thus, one has Heinz
transformation $u$ of the height function $f$ of a minimal graph satisfying%
\[
\frac{1}{\sqrt{1+\left\vert Df\right\vert ^{2}}}\left[
\begin{array}
[c]{cc}%
1+f_{1}^{2} & f_{1}f_{2}\\
f_{2}f_{1} & 1+f_{2}^{2}%
\end{array}
\right]  =D^{2}u\ \ \ \ \ \text{and \ \ }\det D^{2}u=1\ \ \ \text{on
}\mathbb{R}^{2}.
\]

As just obtained, the Hessian $D^{2}u$ of entire solution $u$ is a constant
matrix, and in turn, $Df$ is a constant vector. The original Bernstein theorem
is reached.

In passing, let us note the conjugate function of $f$%
\[
f^{\ast}\left(  x_{1},x_{2}\right)  =\int^{\left(  x_{1},x_{2}\right)  }%
\frac{-f_{2}dx_{1}+f_{1}dx_{2}}{\sqrt{1+\left\vert Df\right\vert ^{2}}}%
\]
satisfies%
\[
Df^{\ast}=\frac{\left(  -f_{2},f_{1}\right)  }{\sqrt{1+\left\vert
Df\right\vert ^{2}}}\ \ \text{and \ \ }\sqrt{1-\left\vert Df^{\ast}\right\vert
^{2}}=\frac{1}{\sqrt{1+\left\vert Df\right\vert ^{2}}}\in\left(  0,1\right)
\]
or%
\[
\frac{Df^{\ast}}{\sqrt{1-\left\vert Df^{\ast}\right\vert ^{2}}}=\left(
-f_{2},f_{1}\right)  .
\]
Consequently%
\[
\operatorname{div}\left(  \frac{Df^{\ast}}{\sqrt{1-\left\vert Df^{\ast
}\right\vert ^{2}}}\right)  =0.
\]

Observe that the above conjugation process from minimal surface to maximal
surface is reversible. \textquotedblleft Incidentally\textquotedblright\ we
have obtained a two dimensional Bernstein type result:\ every entire solution
to the maximal surface equation $\operatorname{div}\left(  Df/\sqrt
{1-\left\vert Df\right\vert ^{2}}\right)  =0$ is linear, which was first
proved up to four dimensions by Calabi, and in general dimension by Cheng-Yau.

\subsubsection{General dimensions}

Next, we outline the argument toward rigidity for semiconvex entire solutions
to $\sigma_{2}\left(  D^{2}u\right)  =1$ in general dimensions.

The Legendre-Lewy transform of a general semiconvex solution satisfies a
uniformly elliptic, saddle equation with bounded Hessian. In the almost convex
case, the new equation becomes concave, thus the Evans-Krylov-Safonov theory
yields the constancy of the bounded new Hessian, and in turn, the old one. To
beat the saddle case, one has to be \textquotedblleft lucky\textquotedblright,
only one time. Recall that, in general the Evans-Krylov-Safonov fails as shown
by the saddle counterexamples of Nadirashvili-Vladuts. Our earlier trace
Jacobi inequality, as an alternative log-convex vehicle, other than the
maximum eigenvalue Jacobi inequality, in deriving the Hessian estimates for
general semiconvex solutions \cite{SY1}, could rescue the saddleness. But the
trace Jacobi only holds for large enough trace of the Hessian. It turns out
that the trace added by a large enough constant satisfies the elusive Jacobi inequality.

Equivalently, the reciprocal of the shifted trace Jacobi quantity is
superharmonic, and it remains so in the new vertical coordinates under the
Legendre-Lewy transformation by a transformation rule. Then the iteration
arguments developed in our joint work with Caffarelli show the
\textquotedblleft vertical\textquotedblright\ solution is close to a
\textquotedblleft harmonic\textquotedblright\ quadratic at one small scale,
\textquotedblleft luckily\textquotedblright\ (two steps in the execution: the
superharmonic quantity concentrates to a constant in measure by applying
Krylov-Safonov's weak Harnack; a variant of the superharmonic quantity, as a
quotient of symmetric Hessian functions of the new potential, is very
pleasantly concave and uniformly elliptic, consequently, closeness to a
\textquotedblleft harmonic\textquotedblright\ quadratic is possible by the
Evans-Krylov-Safonov theory), and the closeness improves increasingly as we
rescale (this is a self-improving feature of elliptic equations, no
concavity/convexity needed). Thus a H\"{o}lder estimate for the bounded
Hessian is realized, and consequently so is the constancy of the new and then
the old Hessian.

Note that, in three dimensions, our proof provides a \textquotedblleft
pure\textquotedblright\ PDE way to establish the rigidity, distinct from the
geometric measure theory way used in our earlier work on the rigidity for
special Lagrangian equations two decades ago.

The details are in the following.

Step 1. Bounded Hessian and uniform ellipticity after Legendre-Lewy transform.

The Legendre-Lewy transform $w\left(  y\right)  =\mathcal{LL}\left[  u\left(
x\right)  +K\left\vert x\right\vert ^{2}/2\right]  $ of a general semiconvex
solution $u\left(  x\right)  $ with $D^{2}u\geq\left(  \delta-K\right)  I$
satisfies a uniformly elliptic, saddle equation with bounded Hessian:%
\begin{gather*}
\left(  x,Du\left(  x\right)  +Kx\right)  =\left(  Dw,y\right)  ,\\
\ 0<D^{2}w=\left(  D^{2}u+K\right)  ^{-1}<\delta^{-1}\ \ \text{or }\lambda
_{i}=\mu_{i}^{-1}-K\geq\delta-K,\\
g\left(  \mu\right)  =-f\left(  \mu^{-1}-K\right)  =-\sigma_{2}\left(
\mu^{-1}-K\right)  =-1.
\end{gather*}
By Lin-Trudinger \cite{LT}, and also Chang-Yuan \cite{ChY}
\[
\lambda_{1}^{-1}\lessapprox f_{\lambda_{1}}\lessapprox\lambda_{1}%
,\ \ \ f_{\lambda_{k\geq2}}\approx\lambda_{1}\ \ \text{for }\lambda_{1}%
\geq\cdots\geq\lambda_{n};
\]
for $\sigma_{2}\left(  \lambda\right)  =1$ with $\lambda_{i}\geq\delta-K,$ all
but one eigenvalues are bounded, $\left\vert \lambda_{k\geq2}\right\vert \leq
C\left(  K\right)  $ and $f_{\lambda_{1}}\lambda_{1}\approx1,$ then%
\[
g_{\mu_{i}}=f_{\lambda_{i}}\mu_{i}^{-2}=f_{\lambda_{i}}\ \left(  \lambda
_{i}+K\right)  ^{2}\approx C\left(  n,K\right)  \ \lambda_{1}\ .
\]
Consequently, level set%
\begin{equation}
\ \left\{  \mu\ |\ g\left(  \mu\right)  =-\sigma_{2}\left(  \mu^{-1}-K\right)
=-1\right\}  \ \text{is\ a uniformly\ elliptic\ surface.} \label{Enue}%
\end{equation}

The new equation also takes the form%
\[
\frac{\sigma_{n-2}\left(  \mu\right)  }{\sigma_{n}\left(  \mu\right)
}-\left(  n-1\right)  K\frac{\sigma_{n-1}\left(  \mu\right)  }{\sigma
_{n}\left(  \mu\right)  }+\frac{n\left(  n-1\right)  }{2}K^{2}=1.
\]

\textbf{Remark.} In the almost convex case $K=\sqrt{2/n\left(  n-1\right)  },$
the new equation becomes%
\begin{equation}
\ \frac{\sigma_{n-1}\left(  \mu\right)  }{\sigma_{n-2}\left(  \mu\right)
}=\left[  \left(  n-1\right)  K\right]  ^{-1}\ , \label{Enc}%
\end{equation}
thus concave. Then the Evans-Krylov-Safonov theory yields the constancy of the
bounded new Hessian%
\[
\left[  D^{2}w\right]  _{C^{\alpha}\left(  B_{R}\right)  }\leq\frac{C\left(
n\right)  }{R^{\alpha}}\left\Vert D^{2}w\right\Vert _{L^{\infty}\left(
B_{2R}\right)  }\leq\frac{C\left(  n\right)  }{R^{\alpha}}\rightarrow
0\ \ \text{as }R\rightarrow\infty,
\]
and in turn, the old one.

The uniformly elliptic level set $\left\{  g\left(  \mu\right)  =-\sigma
_{2}\left(  \mu^{-1}-K\right)  =-1\right\}  $ is saddle for large $K.$ There
are $C^{1,\alpha}$ singular solutions to uniformly elliptic saddle equation in
five dimensions by Nadirashvili-Tkachev-Vladuts.

\bigskip

In the next two steps, we really use the following equivalent uniformly
elliptic saddle equation%
\[
H\left(  D^{2}w\right)  =\sigma_{n}\left(  \mu\right)  \left[  1-\sigma
_{2}\left(  \mu^{-1}-K\right)  \right]  =0,
\]
which in three dimensions, becomes%
\begin{gather*}
-\sigma_{1}\left(  \mu\right)  +2K\sigma_{2}\left(  \mu\right)  -\left(
3K^{2}-1\right)  \sigma_{3}\left(  \mu\right)  =0\ \ \ \text{or}\\
\frac{\sigma_{2}\left(  \mu\right)  }{\sigma_{1}\left(  \mu\right)  }=\left[
2K-\left(  3K^{2}-1\right)  \frac{\sigma_{3}\left(  \mu\right)  }{\sigma
_{2}\left(  \mu\right)  }\right]  ^{-1}.
\end{gather*}

Step 2. Shifted trace Jacobi inequality to rescue saddleness.

For $F\left(  D^{2}u\right)  =\sigma_{2}\left(  \lambda\right)  =1\ $with
$D^{2}u\geq-KI,$ $b\left(  x\right)  =\ln\left(  \bigtriangleup u+nK\right)  $
satisfies the elusive strong subharmonicity%

\[
\bigtriangleup_{F}b=F_{ij}\partial_{ij}b\geq F_{ij}\partial_{i}b\ \partial
_{j}b=\left\vert \nabla_{F}b\right\vert ^{2}%
\]
which is equivalent to the superharmonicity, by a transformation rule in
\cite{SY1}.
\begin{gather*}
\bigtriangleup_{H}a\left(  y\right)  \leq0\\
\text{with\ }a\left(  y\right)  =\frac{1}{\lambda_{1}+K+\cdots+\lambda_{n}%
+K}=\frac{\sigma_{n}\left(  \mu\right)  }{\sigma_{n-1}\left(  \mu\right)
}\overset{n=3}{=}\frac{\sigma_{3}\left(  \mu\right)  }{\sigma_{2}\left(
\mu\right)  }.
\end{gather*}

\textbf{Remark.} For log-convex function $b=\ln\lambda_{\max}\ $or $\ln\left(
\lambda_{\max}+K\right)  ,$ Jacobi inequality $\bigtriangleup_{F}%
b\geq\left\vert \nabla_{F}b\right\vert ^{2}$ holds in three dimensions without
any restriction, in general dimensions with necessary semiconvexity condition.
The reciprocal $\mu_{\min}\left(  y\right)  =\left(  \lambda_{\max}+K\right)
^{-1}$ is superharmonic $\bigtriangleup_{H}\mu_{\min}\left(  y\right)  \leq0.$
But $\mu_{\min}\left(  D^{2}w\right)  $ is not a uniformly elliptic
function/operator on $D^{2}w,$ though concave. Therefore, it is not adequate
via $\mu_{\min}\left(  y\right)  $ to run the Caffarelli-Yuan procedure for
H\"{o}lder of Hessian $D^{2}w.$

For log-linear $b=\ln\bigtriangleup u,$ Qiu showed $\bigtriangleup
_{F}b\overset{\text{3-d}}{\geq}\left\vert \nabla_{F}b\right\vert ^{2}$ in
\cite{Q}. It is indeed another log-convex function%
\[
\ln\bigtriangleup u=\ln\sqrt{\left\vert \lambda\right\vert ^{2}+2}%
\ \ \text{mod }\sigma_{2}\left(  \lambda\right)  =1
\]
satisfying $\bigtriangleup_{F}\ln\bigtriangleup u\geq\left\vert \nabla_{F}%
\ln\bigtriangleup u\right\vert ^{2}$ for large enough $\bigtriangleup u\geq
C\left(  K\right)  $ under semiconvexity $D^{2}u\geq-K.$ It is important to
have the above shifted Jacobi quantity $b\left(  x\right)  =\ln\left(
\bigtriangleup u+nK\right)  $ valid without assuming $\bigtriangleup u$ large
enough, to execute our argument toward constancy of the Hessian.

\textbf{Remark. }The above Jacobi inequality is actually an equality on
minimal surface $\left(  x,f\left(  x\right)  \right)  \in\mathbb{R}^{n}%
\times\mathbb{R}^{1}$ for $\operatorname{div}\left(  Df/\sqrt{1+\left\vert
Df\right\vert ^{2}}\right)  =0:$%
\[
\bigtriangleup_{g}b=\left\vert \nabla_{g}b\right\vert ^{2}+\left\vert
IIg^{-1}\right\vert ^{2}\ \ \text{or \ }\bigtriangleup_{g}\omega
=-\omega\left\vert IIg^{-1}\right\vert ^{2}\leq0,
\]
where $b=\ln\sqrt{1+\left\vert Df\right\vert ^{2}},\ $and$\ \omega
=\left\langle \left(  0,\cdots,0,1\right)  ,N\right\rangle =1/\sqrt
{1+\left\vert Df\right\vert ^{2}}\ $is the effective deformation, while
varying the minimal surface along a Jacobi vector field $J=\left(
0,\cdots,0,1\right)  .$

Step 3. H\"{o}lder estimate of new Hessian on saddle equation.

We illustrate the argument for H\"{o}lder estimates for Hessian on uniformly
elliptic saddle equation in three dimensions, where the idea is not lost, but
the gain is a better understanding of the idea.

Again, the new/\textquotedblleft vertical\textquotedblright\ uniformly
elliptic saddle equation is%
\begin{gather*}
\frac{\sigma_{2}\left(  \mu\right)  }{\sigma_{1}\left(  \mu\right)  }=\left[
2K-\left(  3K^{2}-1\right)  \frac{\sigma_{3}\left(  \mu\right)  }{\sigma
_{2}\left(  \mu\right)  }\right]  ^{-1}\ \ \\
\text{for }1\geq\mu_{i\geq2}\geq c\left(  K\right)  >0\ \ \text{and }1\geq
\mu_{1}>0.
\end{gather*}
Apply Krylov-Safonov's weak Harnack to the bounded superharmonic $a\left(
y\right)  =\frac{\sigma_{3}\left(  \mu\right)  }{\sigma_{2}\left(  \mu\right)
}$ from Step 2, $a\left(  y\right)  $ concentrates to a level $l=\min
_{B_{\varepsilon}}a\left(  y\right)  $ in one small (enough) ball
\ $B_{\varepsilon}\left(  0\right)  ,$ that is $\frac{\sigma_{3}\left(
\mu\right)  }{\sigma_{2}\left(  \mu\right)  }\approx l$ in $99.99999\%$ of
$B_{\varepsilon}.$ Note the concave $\frac{\sigma_{3}\left(  \mu\right)
}{\sigma_{2}\left(  \mu\right)  }$ is not uniformly elliptic, because $\mu
_{1}$ could be close $0.\ $

Step 3 Continued: Remarkably, $\mu$ is approximately on the uniformly elliptic
($1\geq\mu_{3},\mu_{2}\geq c\left(  K\right)  >0$), concave (level set
$\frac{\sigma_{2}\left(  \mu\right)  }{\sigma_{1}\left(  \mu\right)  }=l$ is a
concave surface) equation%
\[
\frac{\sigma_{2}\left(  \mu\right)  }{\sigma_{1}\left(  \mu\right)  }=\left[
2K-\left(  3K^{2}-1\right)  \ l\ \right]  ^{-1}\ \text{\ \ in }%
99.99999\%\ \text{of}\ B_{\varepsilon}%
\]
By existence of Dirichlet problem via Evans-Krylov-Safonov and Alexandrov
maximum principle in measure,
\[
w\left(  y\right)  =\text{quadratic }Q\left(  y\right)  \pm
0.00000000001\ \text{in }B_{\varepsilon/2}.
\]

By self-improving property of the smooth uniformly elliptic equation (no
concavity needed)%
\[
-\sigma_{1}\left(  \mu\right)  +2K\sigma_{2}\left(  \mu\right)  -\left(
3K^{2}-1\right)  \sigma_{3}\left(  \mu\right)  =0,
\]
through iterative procedure $w\left(  y\right)  \approx y^{2+\alpha}$ near
$y=0.$ Similarly near everywhere in $B_{\varepsilon/4}.$ Thus H\"{o}lder
estimates for $D^{2}w.$

Finally, by quadratic scaling%
\[
\left[  D^{2}w\right]  _{C^{\alpha}\left(  B_{R}\right)  }\leq\frac{C\left(
n,K\right)  }{R^{\alpha}}\ \ \overset{R\rightarrow\infty}{\longrightarrow
}\ \ 0.
\]

Thus $D^{2}w$ is constant matrix, in turn, so is $D^{2}u.$

\textbf{Remark.} In the original Caffarelli-Yuan procedure, the superharmonic
quantity $a=\frac{\sigma_{3}\left(  \mu\right)  }{\sigma_{2}\left(
\mu\right)  }$ is $1-e^{K\bigtriangleup u}$ for saddle uniformly elliptic
equation $F\left(  D^{2}u\right)  =0$ with convex level setS along trace $trM$
of the level set $\left\{  F\left(  M\right)  =0\right\}  .$ Once, only one
lucky chance needed, $\bigtriangleup u$ concentrates in measure, solution $u$
is close to a quadratic, by solving the Laplace equation and applying the
Alexandrov maximum principle in measure. Afterwards, the self-improving
machinery of smooth uniformly elliptic equation takes over for H\"{o}lder
estimate of Hessian.

\subsection{Regularity for viscosity solutions}

Any reasonable solution to the $\sigma_{1}$ or Laplace equation
$\bigtriangleup u=1$ such as continuous viscosity in pointwise sense or
distributional solution in integral sense is analytic. The same is not true
for $\sigma_{k\geq3}$ equations, and unknown for $\sigma_{2}$ equation in
dimension five and higher currently. In fact, by now there are
$C^{1,\varepsilon}$ and Lipschitz Pogorelov-like singular viscosity solutions
to $\sigma_{k\geq3}$ equations in dimension three and higher. In the
$\sigma_{n}$ or Monge-Amp\`{e}re equation case, those viscosity solutions are
also singular solutions in the Alexandrov integral sense.

The advance in the joint work with Chen-Shankar \cite{CSY} also led us to
obtain interior regularity (analyticity) for almost convex viscosity solutions
to the quadratic Hessian equation (\ref{Esigma2}), in the joint work with
Shankar \cite{SY2}. Due to similar conceptual and technical challenges--smooth
approximations may not preserve those semiconvexity constraints--we cannot
invoke our available Hessian estimates with Shankar \cite{SY1} for general
semiconvex solutions or with McGonagle and Song \cite{MSY} for almost convex
solutions, while taking the limit and deduce interior regularity. A key
observation is that the Legendre-Lewy transform of any semiconvex viscosity
solution to the equivalent concave equation (\ref{Eec}) stays viscosity
solution to a new concave (\ref{Enc}) (only for the original almost convex
solution) and uniformly elliptic equation (\ref{Enue}) (for all original
semiconvex solutions). In passing, let us note that, for a general fully
nonlinear second order elliptic equation, Alvarez-Lasry-Lions showed that the
Legendre transform of any strictly convex $C^{2}$ solution is a convex
viscosity solution of a conjugate equation. Moreover, the \textquotedblleft
striking\textquotedblright\ role of the $C^{2}$ regularity of the original
solution in their arguments was pointed out \cite[p.281]{ALL}. It follows that
the transformed $C^{1,1}$ solution is smooth by the Evans-Krylov-Safonov
theory. Then the boundedness of the original solutions combined with the
constant rank theorem by Caffarelli-Guan-Ma \cite{CGM} implies that the
original viscosity solution is smooth.

Shortly after, Mooney \cite{M} provided a different proof of the interior
regularity for convex viscosity solutions: every such convex solution is
strictly 2-convex, then all smooth approximated solutions enjoy uniform
Pogorelov-type $C^{1,1}$ and higher derivative estimates by Chou-Wang
\cite{CW}, in turn, the interior regularity by taking limit.

In two dimensions, the above regularity result (now the convexity condition is
automatic) actually also follows from Heinz's famous Hessian estimate earlier.
In three dimensions, the regularity for continuous viscosity solutions to
(\ref{Esigma2}) follows from the Hessian estimate in our joint work with
Warren \cite{WY} and smooth existence with smooth boundary value by
Caffarelli-Nirenberg-Spruck, and also Trudinger. Our most recent joint work
with Shankar \cite{SY4} on Hessian estimates for the quadratic Hessian
equation (\ref{Esigma2}) in four dimensions yields up the same regularity in
four dimensions. There, a direct way to interior regularity without first
deriving the Hessian estimates is also provided. Consequently, a compactness
argument leads to an implicit Hessian estimate.

\subsection{A priori Hessian estimates}

In our long investigation, culminating in the most recent joint work with
Shankar \cite{SY4}, we obtained an implicit Hessian estimate and interior
regularity (analyticity) for the quadratic Hessian equation $\sigma_{2}\left(
D^{2}u\right)  =1$ in four dimensions. Our compactness method (almost Jacobi
inequality--doubling--twice differentiability--small perturbation) also
provides respectively a Hessian estimate for smooth solutions satisfying a
dynamic semiconvexity condition in higher dimensions, which includes
convexity, almost convexity, and semiconvexity conditions appeared in the
recent papers on Hessian estimates, and a non-minimal surface proof for the
corresponding three dimensional results in our earlier joint work with Warren
\cite{WY}.

Other consequence is a rigidity result for entire solutions to the $\sigma
_{2}$ equation with quadratic growth, namely all such solutions must be
quadratic, provided the smooth solutions in dimension $n\geq5$ also satisfying
the dynamic semiconvex assumption.

Again, the Hessian estimate for the $\sigma_{2}$ or Monge-Amp\`{e}re equation
in dimension two was achieved by Heinz in the 1950s. Hessian estimates fail
for the Monge-Amp\`{e}re equation in dimension three and higher, as
illustrated by the famous counterexamples of Pogorelov in the 1970s; those
irregular solutions also serve as counterexamples for cubic and higher order
symmetric $\sigma_{k\geq3}$ equations, for example by Urbas. Hessian estimates
for solutions with certain strict convexity constraints to the
Monge-Amp\`{e}re and $\sigma_{k\geq2}$ equations were derived by Pogorelov and
later Chou-Wang respectively using the Pogorelov technique; some (pointwise)
Hessian estimates in terms of certain integrals of the Hessian were obtained
by Urbas in the early 2000s. The gradient estimates for $\sigma_{k}$ equations
were derived by Trudinger, Chou-Wang in the mid 1990s.

The compactness proof toward an implicit Hessian estimates for almost convex
solutions in \cite{MSY} is based on the concavity of uniformly elliptic
equation (\ref{Enc}) under the Legendre-Lewy transformation, a constant rank
theorem by Caffarelli-Guan-Ma [CGM], on the vertical side, and a strip
argument on the horizontal side.

The proof toward an explicit Hessian estimate for semiconvex solutions is
based on an elusive-Jacobi inequality-satisfying quantity, the maximum
eigenvalue of the Hessian of the solutions,\ envisioned to be true in 2012.
Another essential new device is a mean value inequality for the strongly
subharmonic maximum eigenvalue under the Legendre-Lewy transformation with
uniformly elliptic equation (\ref{Enue}), and its weighted version converted
back to the original variables or horizontal side.

The new idea for Hessian estimates, under a dynamic semiconvexity condition in
dimension five and higher, and consequently interior regularity in dimension
four, is first to get a doubling, or a \textquotedblleft
three-sphere\textquotedblright\ inequality for the Hessian bound on the middle
ball, in terms of Hessian bound on a small inner ball and gradient bound on
the outer large ball:%
\[
\max_{B_{2}\left(  0\right)  }\bigtriangleup u\leq C\left(  r,\left\Vert
u\right\Vert _{Lip\left(  B_{3}\left(  0\right)  \right)  }\right)
\max_{B_{r}\left(  0\right)  }\bigtriangleup u.
\]
Using a Jacobi inequality, true with a lower $\sigma_{3}$-bound condition for
Hessian, satisfied by convex solutions, Guan-Qiu \cite{GQ} reached their
Hessian estimate for the quadratic Hessian equation. Qiu \cite{Q} followed
with his doubling in three dimensions, where the Jacobi inequality was
unconditionally available since \cite{WY}. But the maximum of Guan-Qiu test function%

\begin{gather*}
\left.
\begin{array}
[c]{c}%
P=2\ln\left(  9-\left\vert x\right\vert ^{2}\right)  +\alpha\left\vert
Du\right\vert ^{2}/2+\beta\ \left(  x\cdot Du-u\right) \\
+\ln\max\left\{  \ln\frac{\bigtriangleup u}{\max_{B_{1}\left(  0\right)
}\bigtriangleup u},\gamma^{-1}\right\}
\end{array}
\right. \\
\left.
\begin{array}
[c]{c}%
\text{with\ small\ }\gamma=\gamma\left(  n\right)  >0,\ \text{smaller\ }%
\beta=\beta\left(  \gamma,\left\Vert u\right\Vert _{Lip\left(  B_{3}\left(
0\right)  \right)  }\right)  >0\\
\text{and\ smallest\ }\alpha=\alpha\left(  \gamma,\left\Vert u\right\Vert
_{Lip\left(  B_{3}\left(  0\right)  \right)  }\right)  >0,
\end{array}
\right.
\end{gather*}
could not be ruled out from happening on the small inner ball without the
$\sigma_{3}$-lower bound assumption, thus Qiu's \textquotedblleft
three-sphere\textquotedblright\ inequality.

Now only an almost Jacobi inequality is available in dimension four. In fact,
as observed in 2012, there is no Jacobi inequality in dimension four, and
worse, even no subharmonicity of the Laplace of the log of Hessian in
dimension five, thus the added dynamic semiconvexity condition in higher
dimensions for an almost Jacobi inequality:%
\begin{gather*}
\bigtriangleup_{F}b\overset{n=4}{\geq}\left(  \frac{1}{2}+\frac{\lambda_{\min
}}{\bigtriangleup u}\right)  \left\vert \nabla_{F}b\right\vert ^{2}\geq0;\\
\left.
\begin{array}
[c]{c}%
\bigtriangleup_{F}b\overset{n\geq5}{\geq}\left(  c_{n}+\frac{\lambda_{\min}%
}{\bigtriangleup u}\right)  \left\vert \nabla_{F}b\right\vert ^{2}%
\geq0,\ \ \ \text{\textbf{IF}}\mathbf{\ \ \ }c_{n}+\frac{\lambda_{\min}%
}{\bigtriangleup u}\geq0\ \\
\ \ \ \text{with \ \ \ }c_{n}=\frac{\sqrt{3n^{2}+1}-n+1}{2n}\ \ \text{and
\ \ }b=\ln\bigtriangleup u
\end{array}
\right.  \ \ \text{\ }.
\end{gather*}
Note that for $\sigma_{2}\left(  \lambda\right)  =1,$ we have $\frac
{\lambda_{\min}}{\bigtriangleup u}>-\frac{n-2}{n},$ and at extreme
configuration $\lambda=\left(  K,\cdots,K,-\frac{n-2}{2}K+\frac{1}{\left(
n-1\right)  K}\right)  ,$ one has $\frac{\lambda_{\min}}{\bigtriangleup
u}\underset{\rightarrow}{>}-\frac{n-2}{n}.$ In fact Jacobi inequality holds
$\bigtriangleup_{F}b\overset{n=3}{\geq}\frac{1}{3}\left\vert \nabla
_{F}b\right\vert ^{2}\geq\left(  \frac{1}{2}-\frac{1}{3}\right)  \left\vert
\nabla_{F}b\right\vert ^{2}$ unconditionally in three dimensions.

But the almost Jacobi inequality is really a regular one away from the extreme
configuration of the equation, where the equation is conformally uniformly
elliptic. Qiu's doubling argument can be pushed through.

Now to find a small inner ball where the Hessian is bounded, we first show the
almost everywhere twice differentiability of continuous viscosity solutions,%
\[
u\left(  x\right)  -Q_{y}\left(  x\right)  =o\left(  \left\vert x-y\right\vert
^{2}\right)
\]
by adapting Chauder-Trudinger's argument \cite{CT} for k-convex functions with
$k>n/2,$ with the gradient estimates $\left\Vert Du\right\Vert _{L^{\infty
}\left(  B_{1}\right)  }\leq C\left(  n\right)  \left\Vert u\right\Vert
_{L^{\infty}\left(  B_{2}\right)  },$ actually its integral form of control (a
H\"{o}lder substitute for $k$-convex function was used in \cite{CT}) for
$\sigma_{k}$ equations in \cite{T} and also \cite{CW}, and the fact that
$D^{2}u$ is a bounded Borel measure for solutions of $\sigma_{2}\left(
D^{2}u\right)  =1,$ as%
\[
\int_{B_{1}}\left\vert D^{2}u\right\vert dx\overset{\bigtriangleup
u=\sqrt{2+\left\vert D^{2}u\right\vert ^{2}}}{<}\int_{B_{1}}\bigtriangleup
udx\leq C\left(  n\right)  \left\Vert Du\right\Vert _{L^{\infty}\left(
B_{1}\right)  }.
\]

Then Savin's small perturbation (from the quadratic polynomial at a twice
differentiable point) \cite{S} guarantees the small inner ball with bounded Hessian.

In our most recent follow-up paper with Shankar \cite{SY5}, a new proof of
regularity for strictly convex solutions to $\det D^{2}u=1$ is found, using
similar doubling methods, instead of Euclid distance, now in terms of an
extrinsic distance on the maximal Lagrangian submanifold determined by the
potential Monge-Amp\`{e}re equation. This \textquotedblleft strict
convex\textquotedblright\ regularity was achieved originally by Pogorelov in
the 1960s and 1970s, and generalized by Urbas and Caffarelli in the late 1980s.

\section{Problems}

\noindent\textbf{Problem 1.} Are there singular (Lipschitz) viscosity
solutions, $W^{2,1}$ regularity, and any better partial regularity $\sigma
_{2}\left(  D^{2}u\right)  =1$ in dimension five or higher?

Given the Jacobi inequality is exhausted in our argument for Hessian estimates
and regularity in four dimensions, it is time to look for singular viscosity
solutions and better partial regularity for possible singular viscosity
solutions in dimension five or higher. For example, a dimension estimate on
the singular set of possible singular viscosity solutions. Note that by the
gradient estimate, and then smooth approximations in Lipschitz norm, all
continuous viscosity solutions are Lipschtiz, and by our almost everywhere
twice differentiability \cite{SY4} and Savin's small perturbation \cite{S},
the possible singular set is closed and with zero Lebesgue measure (Also true
for viscosity solutions to $\sigma_{k\geq3}$ equation in dim $n\geq3$).

From $\int_{B_{l}}\left\vert D^{2}u\right\vert dx\overset{\sigma_{2}=1}{<}%
\int_{B_{l}}\bigtriangleup udx$ $\leq\left\Vert Du\right\Vert _{L^{\infty
}\left(  \partial B_{l}\right)  }\left\vert \partial B_{l}\right\vert $ and
the gradient estimate, $D^{2}u$ is a bounded Borel measure. It is reasonable
to expect a $W^{2,1}$ regularity in dim $n\geq5.$ Recall that, unlike function
$\left\vert x_{1}\right\vert ,$ all (convex) viscosity solutions to
$\sigma_{n}\left(  D^{2}u\right)  =1$ have been shown to be $W^{2,1}.$

\noindent\textbf{Problem 2.} Regularity for semiconvex viscosity solutions to
$\sigma_{2}\left(  D^{2}u\right)  =1$ in dimension five or higher.

It is still unclear to us whether semiconvex viscosity solutions are regular,
if only $D^{2}u\geq-KI$ for large $K>0.$ The Legendre-Lewy transform is still
a $C^{1,1}$ viscosity solution of a new uniformly elliptic equation
(\ref{Enue}), for any semiconvex viscosity solution. However, as the new
equation no longer has convex level set, for large $K,$ we are unable to
deduce smoothness for the transformed solution at this point. Without the
smoothness, we are currently unable to obtain a $C^{1,1}$ version of the
constant rank theorem to gain positive definiteness of the semi-positive
Hessian, for the $C^{1,1}$ solution of a uniformly elliptic and inversely
convex equation on the vertical side. Otherwise, the interior regularity for
such semiconvex viscosity solutions would be justified. At this point, it
appears a far stretch to reach regularity for dymanic semiconvex viscosity
solutions in dimension five or higher.

One follow-up of our Hessian estimates for three and the very recent four
dimension $\sigma_{2}$ equation would be

\noindent\textbf{Problem 3.} Derive Schauder and Calder\'{o}n-Zygmund
estimates for variable-right-hand-side equation $\sigma_{2}\left(
D^{2}u\right)  =f\left(  x\right)  $ in dimension four.

With a $C^{1,1}$ assumption on $f\left(  x\right)  ,$ Qiu \cite{Q} has
generalized the arguments in \cite{WY} to obtain Hessian estimates in
dimension three. With an almost sharp Lipschitz assumption on $f\left(
x\right)  ,$ very recently, Zhou \cite{Z} reached the Hessian estimate in
three dimensions, along his Hessian estimates for the \textquotedblleft
twist\textquotedblright\ special Lagrangian equation $\sum_{i=1}^{n}%
\arctan\lambda_{i}/f\left(  x\right)  =c\in\lbrack\left(  n-2\right)
\pi/2,n\pi/2).$ Consequently $C^{2,\alpha}$ estimates follow. Under a small
enough H\"{o}lder seminorm assumption on $f,$ Xu \cite{X} derived interior
$C^{2,\alpha}$ estimates in dimension three. Notice that the interior gradient
estimates in \cite{T} and \cite{CW} needs Lipschitz assumption on $f.$ The
subtle small seminorm constraint is due to the non-uniform elliptic nature of
the equation. The method works in general dimensions such as the recent four
dimension, as long as the interior Hessian estimate is available for the
quadratic Hessian equation with constant right-hand side.

\noindent\textbf{Problem 4. }Any \textquotedblleft
elementary\textquotedblright\ pointwise argument toward Hessian estimates for
$\sigma_{2}\left(  D^{2}u\right)  =1$ in dimension three and four, in general
higher dimension with the dynamic semiconvexity condition, as in the two
dimensional case by Chen-Han-Ou \cite{CHO}, and convex case by Guan-Qiu
\cite{GQ}?

Forthermore, any explicit Hessian bound in terms of the gradient, as the
quadratic exponential dependence in dimension three and general semiconvex case?

To gain more understanding of $\sigma_{2}$ equation, one distinct double
divergence structure of $\sigma_{2}$ from $\sigma_{k\geq3}$ is worth studying.

\noindent\textbf{Problem 5.} Under what additional condition on $u\in
W^{1,2}\left(  \Omega\right)  $ does the equation $\sigma_{2}\left(
D^{2}u\right)  =1$ in the very weak sense (Iwaniec \cite{I}):
\[
\int_{\Omega}\sum_{1\leq i<j\leq n}\left(  \varphi_{ij}u_{i}u_{j}-\frac{1}%
{2}\varphi_{ii}u_{j}^{2}-\frac{1}{2}\varphi_{jj}u_{i}^{2}\right)
dx=\int_{\Omega}\varphi dx\ \ \text{for all\ }\varphi\in C_{0}^{\infty}\left(
\Omega\right)
\]
become$\ \sigma_{2}\left(  D^{2}u\right)  =1,$ say, and $\bigtriangleup u>0,$
in the viscosity sense?

(The double divergence structure is readily seen from the well-known Gauss
curvature formula for graph $\left(  x_{1},x_{2},u\left(  x\right)  \right)
\subset\mathbb{R}^{3}$ with induced metric $g:$ $K=\left(  -\frac{1}%
{2}\partial_{11}g_{22}+\partial_{12}g_{12}-\frac{1}{2}\partial_{22}%
g_{11}\right)  /\left(  \det g\right)  ^{2}=\det D^{2}u/\left(  1+\left\vert
Du\right\vert ^{2}\right)  ^{2}.$)

In dimension two, a (necessary) convexity condition should suffice, also for
$\sigma_{2}\left(  D^{2}u\right)  =1$ in the equivalent Alexandrov sense. In
general dimensions, what about $\bigtriangleup u>0$ in distribution sense for
$u\in C^{1,2^{+}/3}?$ The answer is yes in two dimensions, as shown by Pakzad
\cite{P}. Moreover, in two dimensions, no better than $C^{1,1/3}$
\textquotedblleft very weak\textquotedblright\ solution with sign changing
$\bigtriangleup u$ have been constructed; see the work of Lewicka-Pakzad
\cite{LP} and \cite{CS} \cite{CHI}. It is worth noting that the singular
solution to $\sigma_{2}\left(  D^{2}u\right)  =1$ constructed by C.Y. Li
\cite{L}, $u\left(  x\right)  =\left(  x_{1}^{2}+\cdots+x_{7}^{2}\right)
x_{8}^{7/5}-\frac{25}{84}x_{8}^{3/5}-\frac{25}{28}x_{8}^{14/5}\in
W_{loc}^{1,2}\cap C^{3/5}$ jumps branches, because $\bigtriangleup
u\approx-x_{8}^{-7/5}=\pm\infty$ near $x_{8}=0.$

\end{document}